\documentclass[a4paper]{amsart}

\usepackage{epic, eepic, amsfonts, latexsym, amssymb, graphicx,
multicol, mathrsfs, color, amscd, verbatim, paralist,
xspace, url, euscript, stmaryrd,  amsmath, enumitem,
bbold, multirow, tikz, mathtools}
\usepackage[all,pdf,cmtip]{xy}

\usepackage[colorlinks, linkcolor=blue, citecolor=magenta, urlcolor=cyan]{hyperref}

\usepackage{graphicx}

\usepackage{float}

\def\hat{\widehat}

\def\RR{{\mathbb R}}
\def\CC{{\mathbb C}}

\def\NN{{\mathbb N}}

\newcommand{\suchthat}{\;\ifnum\currentgrouptype=16 \middle\fi|\;}

\def\ra{\rightarrow}

\def\hat{\widehat}

\def\qed{{\hfill $\Box$}}
\newtheorem{theorem}{THEOREM}[section]

\theoremstyle{definition}

\theoremstyle{remark}
\newtheorem{remark}[theorem]{Remark}

\makeatletter
\def\blfootnote{\xdef\@thefnmark{}\@footnotetext}
\makeatother

\begin{document}

\title[Hyperbolicity of tube domains in $\CC^2$]{On the Kobayashi hyperbolicity
\vspace{0.1cm}\\
of tube domains in $\CC^2$}
\author[Isaev]{Alexander Isaev}

\address{Mathematical Sciences Institute\\
Australian National University\\
Acton, ACT 2601, Australia}
\email{alexander.isaev@anu.edu.au}

\maketitle

\thispagestyle{empty}

\pagestyle{myheadings}

\begin{abstract} We construct elementary counterexamples to the criterion for Kobayashi hyperbolicity for a class of tube domains in $\CC^2$ proposed by J.-J. Loeb.
\end{abstract}

\section{Introduction}\label{intro}
\setcounter{equation}{0}

A connected complex manifold $X$ is said to be {\it Kobayashi-hyperbolic}\, if the Koba\-ya\-shi pseudodistance on $X$ is in fact a distance (see \cite{K} for details). For $X$ endowed with a Riemannian metric, hyperbolicity is equivalent to the following property: for any point $x\in X$ there exist a neighborhood $U$ of $x$ and a constant $M>0$ such that for all holomorphic maps $f:\Delta\ra X$ with $f(0)\in U$ one has $||df(0)||<M$, where $\Delta$ is the unit disk in $\CC$. Verification of hyperbolicity may be a difficult task even for very special classes of manifolds.

In this note we focus on {\it tube domains}\, in $\CC^n$, i.e., domains of the form\linebreak $T_D:=D+i\RR^n$, where $D$ is a domain in $\RR^n$ called the {\it base}\, of $T_D$. As pointed out in \cite{L} (see also Theorem 13.6.2 in \cite{JP}), for a tube domain $T_D\subset\CC^n$ the hyperbolicity property is equivalent to the following condition: for any point $x\in D$ there exist a neighborhood $U$ of $x$ in $D$ and a constant $M>0$ such that for all harmonic maps $f:\Delta\ra D$ with $f(0)\in U$ one has $||df(0)||<M$. 

From now on, we assume that $n=2$. It is somewhat surprising that so far no easily verifiable criterion for the hyperbolicity of a tube domain has been found even in this situation. By Bochner's theorem, the envelope of holomorphy of $T_D$ coincides with $T_{\hat D}$, where $\hat D$ is the convex hull of $D$ (see, e.g., Section 21 in \cite{V}), and it is natural to investigate hyperbolicity separately in each of the cases:\linebreak (i) $T_{\hat D}\not=\CC^2$ and (ii) $T_{\hat D}=\CC^2$. In \cite{HI} we presented several classes of hyperbolic domains in $\CC^2$ falling in case (ii). For example, we showed that $T_D$ is hyperbolic  if $D$ is a domain bounded by two spirals, where a spiral is a curve defined in polar coordinates in $\RR^2$ by the equation $r=g(\varphi)$, with $g$ being an increasing function of $\varphi$ such that $\lim_{\varphi\ra-\infty}g(\varphi)=0$ and $\lim_{\varphi\ra+\infty}g(\varphi)=\infty$. However, there is no comprehensive description of all hyperbolic domains covered by case (ii) (cf.~Question 13.6 in \cite{JP}). On the other hand, for domains in $\CC^2$ falling in case (i) a hyperbolicity criterion was proposed by J.-J. Loeb in \cite{L}, as stated below.

Let $D\subset\RR^2$ be a domain with $\hat D\ne\RR^2$. Writing coordinates in $\CC^2$ as $z_j=x_j+iy_j$, $j=1,2$, we may assume without loss of generality that $D$ lies in the half-space $\{x_2>0\}$. Then the result of \cite{L} asserts that $T_D$ is hyperbolic if and only if there is no point $(a_1,a_2)\in D$ for which there exists a sequence of real numbers $\{b_k\}$ convergent to $a_2$ with the property that the segment $[-k,k]\times\{b_k\}$ lies in $D$ for all $k\in\NN$. The necessity implication is obvious. Regarding the sufficiency implication, as M.~Jarnicky and P.~Pflug observed, the argument provided in \cite{L} only yields the following weaker statement (see part (b) of Theorem 13.6.6 in \cite{JP}):

\begin{theorem}\label{correctedthm} If $T_D$ is not hyperbolic then there exists a point $(a_1,a_2)\in D$ such that for every $k\in\NN$ one can find a real-analytic function $\gamma_k(t)$ on $[-k, k]$, with $(t,\gamma_k(t))\in D$ and $|\gamma_k(t)-a_2|\le 1/k$ for all $t$.
\end{theorem} 

On the other hand, to the best of our knowledge, no counterexample to the sufficiency implication of Loeb's theorem has been found so far (cf.~part (a) of Remark 13.6.7 in \cite{JP}). In this note we construct such a counterexample thus clearing the confusion that has existed around Loeb's result for a  number of years. Namely, we show:

\begin{theorem}\label{result} There exists $D\subset\RR^2$ lying in the half-space $\{x_2>0\}$ such that $T_D$ is not hyperbolic and for $b\in\RR$ and $k\in\NN$ no segment $[-k,k]\times\{b\}$ is contained in $D$. Such a domain $D$ can be chosen to have a $C^{\infty}$-smooth boundary.
\end{theorem}

Although the examples provided below are elementary, they are nevertheless surprising as one does expect that obstructions for the hyperbolicity of tube domains should indeed be in some sense \lq\lq linear\rq\rq\, (cf.~Loeb's statement). The idea behind the examples inspires a partial converse to Theorem \ref{correctedthm}, which further emphasizes the theme of the \lq\lq linearity\rq\rq\, of obstructions:

\begin{theorem}\label{converselinear} Let $D\subset\RR^2$ be a domain lying in the half-space $\{x_2>0\}$ and satisfying the following condition: there exists a point $(a_1,a_2)\in D$ with the property that for every $k\in\NN$ one can find an affine function $\gamma_k(t)=c_kt+d_k$ such that $(t,\gamma_k(t))\in D$ and $|\gamma_k(t)-a_2|\le 1/k$ for all $t\in [-k,k]$. Then $T_D$ is not hyperbolic.
\end{theorem}

However, it is not clear from the proof of Theorem \ref{correctedthm} given in \cite{JP} that one can always choose the function $\gamma_k$ appearing there to be affine. Thus, the problem of eliminating the gap between necessary and sufficient conditions for the hyperbolicity of tube domains in $\CC^2$ whose envelope of holomorphy is not all of $\CC^2$ remains open.

{\bf Acknowledgement.} This work is supported by the Australian Research\linebreak Council.

\section{The examples and Proof of Theorem \ref{converselinear}}\label{sect1}
\setcounter{equation}{0}

We start by construct domains satisfying the requirements of Theorem \ref{result}. First, let
$$
D:=\{(x_1,x_2)\in\RR^2\mid 0<x_2<2\}\setminus\Bigl(\{-1\}\times[1,2]\cup\{1\}\times[0,1]\Bigr)
$$
as shown in Fig. 1 below. Clearly, for $b\in\RR$ and $k\in\NN$ no segment $[-k,k]\times\{b\}$ is contained in $D$. 

\begin{figure}
\begin{center}  
\includegraphics[width=4.2in]{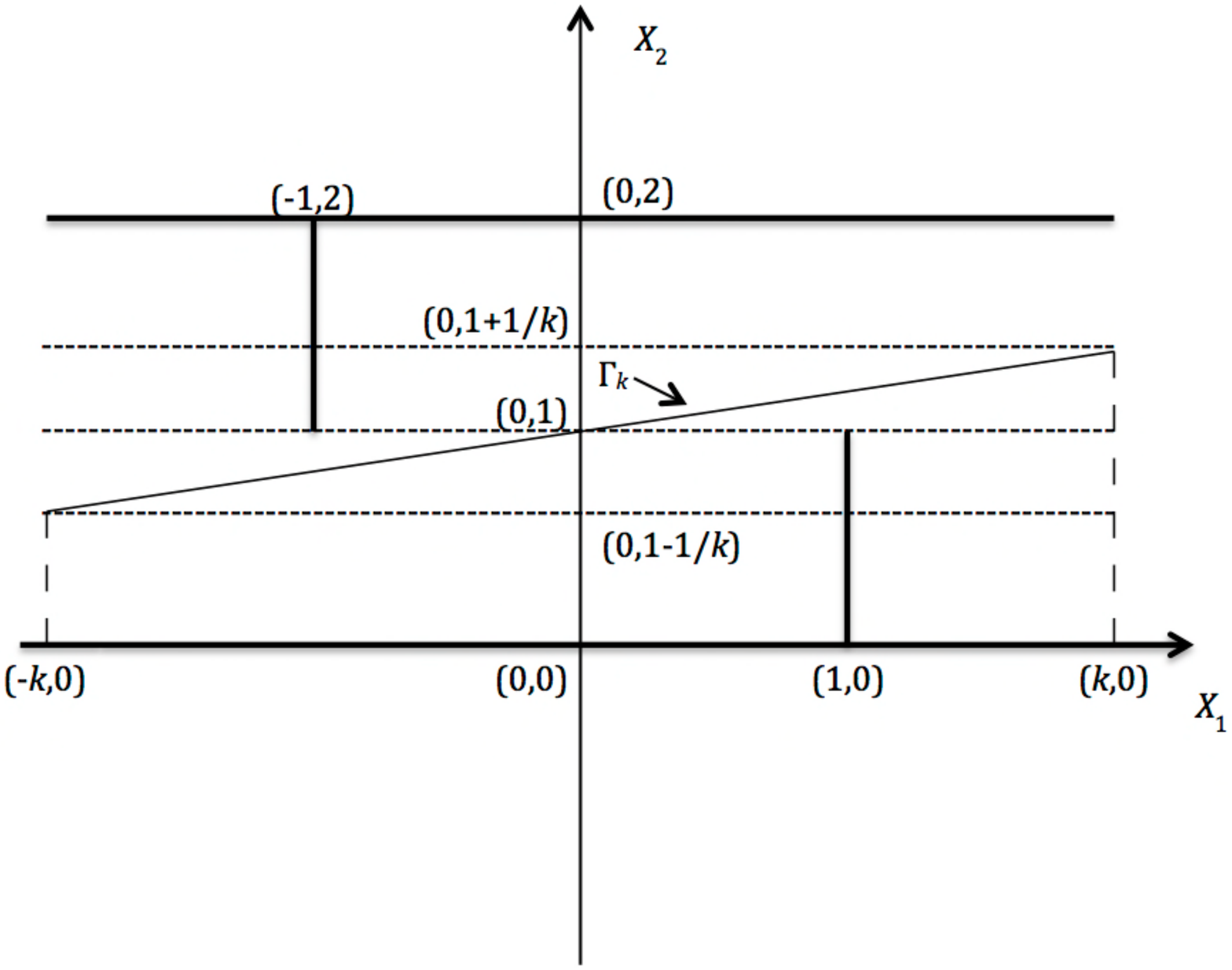}
\caption{\small \sl }  
\end{center}  
\end{figure} 

We will now prove that $T_D$ is not hyperbolic. Let $a:=(0,1)\in D$. We will construct a sequence of holomorphic mappings $f_k:\Delta\to T_D$ such that $f_k(0)=a$ and $||df_k(0)||\ra\infty$ as $k\ra \infty$. Define
$$
f_k:\Delta\to T_D,\quad z\mapsto \left(kz, \frac{1}{k}z+1\right).
$$
Clearly, $f_k(0)=a$ and
$$
df_k(0)=\left(k,\frac{1}{k}\right).
$$
Hence, $||df_k(0)||\ra\infty$ as $k\ra \infty$, which shows that $T_D$ is not hyperbolic.

The above example can be modified by choosing
$$
D:=\{(x_1,x_2)\in\RR^2\mid 0<x_2<2\}\setminus (S_1\cup S_2),
$$
where 
$$
S_1\subset \{(x_1,x_2)\in\RR^2\mid x_1\le 0,\, 1\le x_2\le 2 \}
$$
is a closed region whose boundary contains a curve joining a pair of points on the line $\{x_2=2\}$ and passing through a point on the line $\{x_2=1\}$, and  
$$
S_2\subset \{(x_1,x_2)\in\RR^2\mid x_1\ge 0,\, 0\le x_2\le 1 \}
$$
is a closed region whose boundary contains a curve joining a pair of points on the line $\{x_2=0\}$ and passing through a point on the line $\{x_2=1\}$ as shown in Fig. 2 below. 

\begin{figure}[H]
\begin{center}  
\includegraphics[width=4.2in]{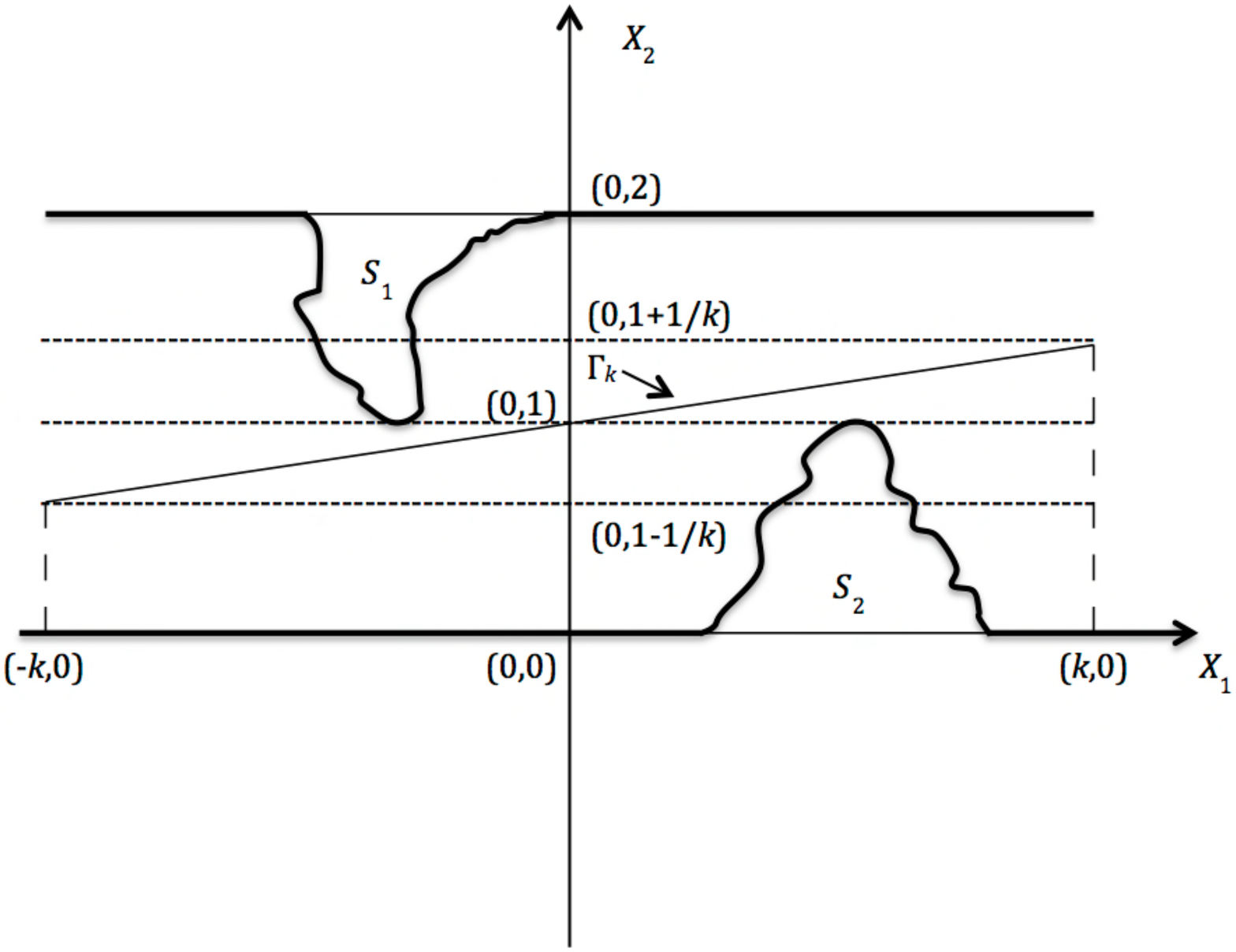}
\caption{\small \sl }  
\end{center}  
\end{figure} 

\noindent Clearly, $S_1$, $S_2$ can be chosen to ensure that $\partial D$ is smooth. Moreover, for any domain of this kind there exists $k_0\in\NN$ such that for $b\in\RR$ and $k\ge k_0$ no segment $[-k,k]\times\{b\}$ is contained in $D$. It is easy to make a choice of $S_1$, $S_2$ so that $k_0=1$, which completes the proof of Theorem \ref{result}.\qed

\begin{remark}\label{gamma} Observe that the examples given above satisfy the condition stated in Theorem \ref{correctedthm} with $a=(0,1)$ and
$$
\gamma_k(t)=\frac{1}{k^2}t+1.
$$
The curve $\Gamma_k(t):=(t,\gamma_k(t))$, with $t\in [-k,k]$, is a line segment as shown in Figs. 1, 2.
\end{remark}

The proof of Theorem \ref{converselinear} is based on a similar idea. We will construct a sequence of holomorphic mappings $f_k:\Delta\to T_D$ (with $k$ sufficiently large) such that\linebreak $f_k(0)\ra a$ and $||df_k(0)||\ra\infty$ as $k\ra \infty$. Define
$$
f_k:\Delta\to \CC^2,\quad z\mapsto \left(a_1+\frac{k}{2}z, c_k\left(a_1+\frac{k}{2}z\right)+d_k\right).
$$
Clearly, $f_k(\Delta)$ lies in $T_D$ if $k$ is large enough. Now observe that $c_k\to 0$ and $d_k\to a_2$, which yields $f_k(0)\to a$ as $k\to\infty$. Furthermore,
$$
df_k(0)=\frac{1}{2}\left(k,kc_k\right),
$$
hence $||df_k(0)||\ra\infty$ as $k\ra \infty$. This shows that $T_D$ is not hyperbolic as\linebreak required. \qed

\end{document}